\documentclass{amsart}%
\usepackage{amssymb}
\usepackage{amsfonts}
\usepackage{amsmath}
\usepackage{graphicx}%
\newtheorem{theorem}{Theorem}
\theoremstyle{plain}

\newtheorem{definition}{Definition}
\newtheorem{example}{Example}

\newtheorem{lemma}{Lemma}

\numberwithin{equation}{section}
\begin{document}
\title[Range of a StOp]{The Range of a Steiner Operation}
\author{L. H. Harper}
\address{Department of Mathematics\\
University of California, Riverside\\
Riverside, CA 92521}
\email{harper@math.ucr.edu}
\date{March 8, 2012}

\begin{abstract}
This paper answers a fundamental question in the theory of Steiner operations
(StOps) as defined and studied in the monograph, \cite{Har04}. StOps are
morphisms for combinatorial isoperimetric problems, analogous to Steiner
symmetrization for continuous isoperimetric problems. The usefulness of a
StOp, $\varphi:\mathbf{2}^{V}\rightarrow\mathbf{2}^{V}$, $V$ a finite set,
depends on having an efficient representation of its range. In \cite{Har04}
the problem was treated case-by-case. In each case the StOp induced a partial
order, $\mathcal{P}$, on $V$ so that $Range\left(  \varphi\right)
=\mathcal{I}\left(  \mathcal{P}\right)  $, the set of all order ideals of
$\mathcal{P}$. Here we show (directly from the axioms for a StOp) that every
idempotent StOp admits such a representation of its range ($\mathcal{P}$ is
then called the StOp-order of $\varphi$). That result leads to another
question: What additional structure does $Range\left(  \varphi\right)  $ have?
The answer is none. We show that every finite poset is the StOp-order of some
idempotent Steiner operation.

\end{abstract}
\maketitle

\section{Background}

\subsection{Combinatorial Isoperimetric Problems}

Certain combinatorial optimization problems are analogous to the classical
isoperimetric problem of plane geometry. The simplest example is the
\textit{edge-isoperimetric problem} (EIP) on a graph, $G=\left(  V;E\right)
$. $V$ is the (finite) set of \textit{vertices} of $G$\textit{\ }and
$E\subseteq\binom{V}{2}$, (unordered pairs of (distinct) vertices) is its set
of \textit{edges}. The \textit{edge-boundary} of a set $S\subseteq V$ is
$\Theta\left(  S\right)  =\left\{  e\in E:e=\left\{  u,v\right\}  ,u\in
S,v\notin S\right\}  $, the edges having one end in $S$ and the other end in
$V-S$. The EIP on $G$ is to minimize $\left\vert \Theta\left(  S\right)
\right\vert $ given $\left\vert S\right\vert $. In general the EIP is a hard
problem (NP-complete), but certain special cases of interest for computer
science and engineering have been solved. The author's first paper, written
fifty years ago (1962), solved the EIP for the graph of the $n$-cube.

\subsection{Steiner Operations}

\subsubsection{Axioms}

The concept of Steiner operation (\textit{StOp}) was the result of a
deliberate effort, beginning around 1976, to identify morphisms for
combinatorial isoperimetric problems (see \cite{Har77}). The basic theory and
results that followed were surveyed in \cite{Har04}, published in 2004. There,
on pages 27-28 a StOp is characterized as a set-map, $\varphi:\boldsymbol{2}%
^{V}\rightarrow\boldsymbol{2}^{V}$, $V$ being a finite set with
\textit{boundary function}, $\Omega$, having the following properties:

\begin{enumerate}
\item $\varphi$ preserves the size of subsets: $\forall S\subseteq V$,
$\left\vert \varphi\left(  S\right)  \right\vert =\left\vert S\right\vert $,

\item $\varphi$ does not increase the size of their boundaries: $\left\vert
\Omega\left(  \varphi\left(  S\right)  \right)  \right\vert \leq\left\vert
\Omega\left(  S\right)  \right\vert $,

\item $\varphi$ preserves the structure of $\boldsymbol{2}^{V}$: $\forall
S\subseteq T\subseteq V$, $\varphi\left(  S\right)  \subseteq\varphi\left(
T\right)  $\textit{.}
\end{enumerate}

The original family of StOps \cite{Har77}, called \textit{stabilization}, was
derived from reflective symmetry of a graph embedded in Euclidean space,
$\mathbb{R}^{n}$ (\textit{e.g. }the graph of the $n$-cube). G.-C. Rota's
observation that stabilization is a discrete analog of Steiner symmetrization
(See \cite{P-S} and \cite{B-Z} for background on Steiner symmetrization or
Google "Steiner symmetrization") added considerable \textit{gravitas} to the
project. Later it was discovered that stabilization is also a StOp for the
vertex-isoperimetric problem (VIP) on the graph, $G=\left(  V;E\right)  $ (the
boundary, $\Phi\left(  S\right)  $, is the set of vertices in $V-S$ with
neighbors in $S$). Also there is another whole class of systematic StOps,
called \textit{compression}. Compression is based on a product decomposition
of $G$, one of the factors having nested solutions for the isoperimetric
problem. Unlike stabilization, compression had been discovered independently
many times, appearing in the majority of papers on combinatorial isoperimetric
problems (See \cite{Har04} or \cite{Eng} for more background and details).
Furthermore, Steiner symmetrization and most of its variants (such as Schwartz
symmetrization) are compressions.

Properties (1), (2) \& (3) above may be regarded as axioms for StOps. However,
in the light of experience and expedience, we wish to extend these axioms a
bit. First by allowing the domain and codomain of a StOp to be the more
general $\mathcal{I}\mathfrak{(}\mathcal{P)}$, where $\mathcal{P}=\left(
V;\leq\right)  $, $\leq$ being a partial order relation on $V$ and
$\mathcal{I}\mathfrak{(}\mathcal{P)}$ the set of all ideals of $\mathcal{P}$
(for definitions see Sec. 2.1). In scheduling problems, where $V$ is a set of
tasks, $\mathcal{P}$ represents precedence constraints under which the tasks
must be performed. No task may be worked on until all its antecedents have
been completed. Also, as mentioned in our abstract and demonstrated in
\cite{Har04}, such restrictions on the domain of a StOp arise from Steiner
operations themselves.

StOps were created to help solve combinatorial isoperimetric problems. These
problems may, in principle, be solved by Brute Force, trying out all
possibilities in the domain, $\mathcal{I}\mathfrak{(}\mathcal{P)}$. However,
if $\mathcal{I}\mathfrak{(}\mathcal{P)}$ is large, the cost may be
prohibitive. The range of $\varphi$ will generally be much smaller than its
domain, but to be able to take advantage of the reduction in size we must have
a simple way to identify sets in
\[
Range\left(  \varphi\right)  =\left\{  T\in\mathcal{I}\mathfrak{(}%
\mathcal{P)}:\exists S\in\mathcal{I}\mathfrak{(}\mathcal{P)}\text{ such that
}\varphi\left(  S\right)  =T\right\}  \text{,}%
\]
and to generate them all efficiently. This fundamental technical problem in
the theory of StOps was apparent right from the start, with stabilization, and
an effective answer was found (see \cite{Har77}). For the edge- and
vertex-isoperimetric problems on graphs, the initial domain is $\boldsymbol{2}%
^{V}$, all subsets, so $\mathcal{P}$ $=\Delta,$ the discrete order on $V$. It
was observed that each of the basic stabilization operations, defined by a
reflective symmetry and a point not on the fixed hyperplane of the reflection
(the Fricke-Klein point), induced a simple partial order on $V$. The basic
stabilizations, with a common Fricke-Klein point so that they all had a common
total extension (are \textit{consistent}), could be composed, giving a StOp
that combined the simplifications of its constituents. Surprisingly, although
basic stabilizations were idempotent (\textit{i.e. }$\varphi^{2}=\varphi)$,
their compositions were generally not. However, if repeatedly (cyclically)
composed, they would eventually become constant and therefore idempotent. The
range of this superstabilization was then characterized (Theorem 4 of
\cite{Har77}) as $\mathcal{I}\mathfrak{(}\mathcal{Q)}$, $\mathcal{Q}$ being
the transitive closure of the union of all the basic partial orders.

Later on the author was pleasantly surprised to note that compression fit into
the same theoretical framework: Basic compressions are idempotent Steiner
operations, each defining a partial order on $V$. If consistent, compressions
could be cyclically composed to give an idempotent StOp combining all their
simplifications into one. The range of that supercompression is exactly the
set of all ideals of the transitive closure of the union of the basic partial
orders (see \cite{Har04}, Section 3.3.4). Also, all \textit{ad hoc}
(unsystematic) Steiner operations have been found to have partial orders
representing their ranges.

All of these StOps and their associated StOp-orders share the property that
$\varphi\left(  S\right)  $ is "lower" (wrt the StOp-order) than $S$. In each
case this is justified by a one-to-one function, $f_{S}:S\rightarrow
\varphi\left(  S\right)  $ such that $\forall x\in S,f_{S}\left(  x\right)
\leq x$, the definition of $f_{S}$ depending on the definition of
$\varphi\left(  S\right)  $. This leads us to add a fourth axiom:

\begin{enumerate}
\item[(4)] $\exists$a total extension
\[
\tau:\mathcal{P}\rightarrow\mathbf{n=}\left\{  0<1<...<n-1\right\}  ,
\]
($\tau$ one-to-one and onto so $\left\vert V\right\vert =n$ and $x\leq
y\Rightarrow\tau\left(  x\right)  \leq\tau\left(  y\right)  $). Also, with
$\tau\left(  S\right)  $ defined to be $%
{\displaystyle\sum\limits_{x\in S}}
$ $\tau\left(  x\right)  $, $\forall S\in\mathcal{I}\mathfrak{(}\mathcal{P)},$
$\tau\left(  \varphi\left(  S\right)  \right)  \leq\tau\left(  S\right)  $ and
$\tau\left(  \varphi\left(  S\right)  \right)  =\tau\left(  S\right)
\Rightarrow\varphi\left(  S\right)  =S${\LARGE .}
\end{enumerate}

For stabilizations, $\tau$ is the order induced by proximity to the
Fricke-Klein point. For compressions it is the total order of the inductive
hypothesis. If a set of StOps share a common $\tau$ they are called
\textit{consistent}. Axiom 4 makes the repeated compositions of consistent
StOps eventually constant since $\tau\left(  \varphi\left(  S\right)  \right)
$ can only decrease a finite number of times. The requirement of consistency
is a pragmatic one. Compositions of StOps that are not consistent would still
be StOps but their composition could cycle and may not become constant. The
reduction in size of the range achieved by the composition of consistent StOps
is seems to be more than that achieved by identifying equivalence classes of
sets. A good example of that is the solution of the edge-isoperimetric problem
on $V_{600}$, the graph of the 600-vertex regular solid in 4 dimensions (See
\cite{Har04}). The Brute Force solution would generate all $2^{600}%
\gtrsim10^{180}$subsets of vertices, an impossible task. The symmetry group of
$V_{600}$ is of order $14,400$ (See Coxeter's classic monograph \cite{Cox}) so
there are at least $10^{180}/14,400\gtrsim\allowbreak6.9\times10^{177}$
equivalence classes of sets of vertices. That is still a huge number and it is
not even clear how to generate those equivalence classes efficiently. However,
$V_{600}$ has $60$ reflective symmetries and the superstabilization they
generate has about $\ 10^{10}$ sets in its range. Those sets are ideals in the
stabilization-order (aka the Bruhat order) of $V_{600}$ and can be recursively
generated in lexicographic order very efficiently. A Brute Force solution of
the edge-isoperimetric problem on $V_{600}$ (generating all $10^{10}$ sets in
the range of its superstabilization) was carried out on a 3 MgHz PC in one day.

We also modify the second axiom by extending "size of the boundary" to any
functional, $\partial:\boldsymbol{2}^{V}\rightarrow\mathbb{R}$, requiring that

\begin{enumerate}
\item[(2$%
\acute{}%
$)] $\varphi$ does not increase the size of their boundaries: $\partial\left(
\varphi\left(  S\right)  \right)  \leq\partial\left(  S\right)  $.
\end{enumerate}

The author's monograph \cite{Har04} treats concrete Steiner operations
(StOps), principally stabilization \& compression, and the partial orders
(StOp-orders) that characterize their ranges. These concepts systematically
simplify hard problems, pointing the way to subsequent developments such as
passage to a continuous limit. Rather than survey the whole book here, we just
give a glimpse of the culminating application, the solution of a problem posed
by A. A. Sapozhenko. It demonstrates the power and efficiency of Steiner
operations and their StOp-orders. Sapozhenko asked about the VIP on the
Johnson graph, $J(d,n)$. The vertices of $J(d,n)$ are $n$-tuples of 0s \& 1s
with exactly $d$ 1s. Two such vertices are neighbors if they differ in exactly
two places. $J(d,n)$ does not have nested solutions for $d>1$. $J(d,n)$ is not
factorable as a product, so even the nested solutions of $J(1,n)$ cannot be
used for compression. The most successful strategy for solving combinatorial
isoperimetric problems without nested solutions has been to pass to a
continuous limit and apply calculus (See Chapter 10 of \cite{Har04}). This
works for several other problems lacking nested solutions: The EIP on $\left(
\mathbb{Z}_{n}\right)  ^{d}$ (the $d$-fold product of $n$-cycles) and the VIP
on $\left(  \mathbb{K}_{n}\right)  ^{d}$, (the $d$-fold product of complete
graphs on $n$ vertices). The (1-dimensional) compression-order for both
$\left(  \mathbb{Z}_{n}\right)  ^{d} $ \& $\left(  \mathbb{K}_{n}\right)
^{d}$ is $\boldsymbol{n}^{d}$ and the limit of $\left(  \boldsymbol{n}%
/n\right)  ^{d}$ as $n\rightarrow\infty$ is $\left[  0,1\right]  ^{d}$, the
unit $d$-cube. Of course $\left[  0,1\right]  ^{d} $ has different boundaries
for the limits of the EIP \& VIP. Bollobas \& Leader solved the EIP on
$\left[  0,1\right]  ^{d}$ with a discontinuous modification of compression
(relaxation of Axiom 3 and induction on $d$). The present author solved the
VIP on $\left[  0,1\right]  ^{d}$ with a discontinuous modification of
stabilization. It is not apparent how to pass to a continuous limit with
$J(d,n)$ but stabilization with respect to the symmetric group acting on its
coordinates transforms it so that the limit becomes obvious. The limit of the
stabilization-order of $J(d,n)$ (as $n\rightarrow\infty$) is the continuous
poset,
\[
\mathcal{L}\left(  d\right)  =\left\{  x\in\left[  0,1\right]  ^{d}:x_{1}\leq
x_{2}\leq...\leq x_{d}\right\}  ,
\]
ordered coordinatewise (See Section 10 of \cite{Har04} for details). The
solution of Sapozhenko's problem follows immediately from the symmetry of the
solution of the VIP on $\left[  0,1\right]  ^{d}$. The logic of the solution
is the same as that given in elementary books for the reduction of Dido's
problem to the classical isoperimetric problem in the plane, a "Didonean
embedding" of $\mathcal{L}\left(  d\right)  $ into $\left[  0,1\right]  ^{d}$
(See \cite{Har04}, Chapter 10 for details). One might expect that the EIP on
$J(d,n)$ could be similarly solved. However, the solution of the EIP on
$J(d,n)$ is not symmetric under interchange of coordinates, so the embedding
is not Didonean. Despite considerable effort, the author has not been able to
adapt the techniques that solved the EIP \& VIP on $\left[  0,1\right]  ^{d}$
to solve the EIP on $\mathcal{L}\left(  d\right)  $ (\textit{i.e. }%
$\lim_{n\rightarrow\infty}J(d,n)$). It remains an open problem.

\subsubsection{Is there a theorem here?}

With all those "coincidences", the range of so many different Steiner
operations, $\varphi:\mathcal{I}\mathfrak{(}\mathcal{P)\rightarrow
I}\mathfrak{(}\mathcal{P)}$ being represented as $\mathcal{I}\mathfrak{(}%
\mathcal{Q)}$ for some $\mathcal{Q\supseteq P}$, it was natural to wonder if
\textit{every} StOp has such a partial order characterizing its range? We kept
coming back to this question because of the efficacy and power of StOp-orders,
but each time were brought up short by the lack of any obvious source for the
additional order relations. Where could they possibly come from? Then one day
we realized that there already was a precedent in the literature for just such
spontaneous creation of order: Garret Birkhoff's characterization of finite
distributive lattices. (See Birkhoff's classic monograph, \underline{Lattice
Theory} \cite{Bir}, Section III.3).

\section{The Ideal Transform and its Ramifications}

Birkhoff's theorem requires some background. We now summarize the definitions
and basic results for it. We have taken these from the monograph by Davey \&
Priestley \cite{D-P} to which the reader may refer for proofs and additional
theory. Also see Gratzer's more recent monograph \cite{Gra}.

\subsection{Posets and Ideals}

A \textit{partial order}, $\leq$, on a set, $V$, is a binary relation, $\leq$
$\subseteq V\times V$, which is

\begin{enumerate}
\item Reflexive: $\forall x\in V,x\leq x$,

\item Antisymmetric: $\forall x,y\in V$, $\left(  x\leq y\text{ }\&\text{
}y\leq x\right)  \Rightarrow\left(  x=y\right)  $,

\item Transitive: $\forall x,y,z\in V$, $\left(  x\leq y\text{ }\&\text{
}y\leq z\right)  \Rightarrow\left(  x\leq z\right)  .$
\end{enumerate}

A \textit{partially ordered set} (\textit{poset}), $\mathcal{P}=\left(
V;\leq\right)  $, consists of a set, $V$, with a partial order, $\leq$, on $V$.

\begin{example}
$\boldsymbol{n}=\left\{  0<1<...<n-1\right\}  $ is a total order (chain) of
size $n$.
\end{example}

\begin{example}
$n=\left\{  0,1,...,n-1\right\}  $ is a discrete order (antichain) of size
$n$. In this case the partial order is $\Delta_{n}=\left\{  \left(
0,0\right)  ,\left(  1,1\right)  ,...,\left(  n-1,n-1\right)  \right\}  $, the
identity relation.
\end{example}

\begin{example}
The Boolean lattice, $\mathcal{B}_{n}$, with $n$ generators is $\mathbf{2}%
^{n}=\left\{  0<1\right\}  \times\left\{  0<1\right\}  \times...\times\left\{
0<1\right\}  $, ordered coordinatewise. $\mathcal{B}_{n}$ is isomorphic to the
power set of an $n$-set.
\end{example}

A set, $I\subseteq V$, is called an (order) \textit{ideal} of $\mathcal{P}%
=\left(  V;\leq\right)  $, if $\left(  y\in I\text{ }\&\text{ }x\leq y\right)
\Rightarrow\left(  x\in I\right)  $. In \cite{D-P} these are called down-sets,
in \cite{Bir} hereditary sets.

\begin{example}
$\boldsymbol{m}$ is an ideal of $\boldsymbol{n}$ $\Leftrightarrow m\leq n$.
\end{example}

$\mathcal{I}\mathfrak{(}\mathcal{P)}=\left\{  I\subseteq V:I\text{ is an ideal
of }\mathcal{P}\right\}  $ is called the \textit{ideal-set} or \textit{ideal
transform} of $\mathcal{P}$. In \cite{D-P} our $\mathcal{I}\mathfrak{(}%
\mathcal{P)}$ is denoted $\mathcal{O}\mathfrak{(}\mathcal{P)}$, $\mathcal{I}$
being reserved for the ideal-sets of lattices, which have additional
structure. We use $\mathcal{I}$ for both, feeling that the concept of ideal
for posets, lattices (and rings) are essentially the same, differing only by
context (different categories).

\begin{example}
$\mathcal{I}\mathfrak{(}\boldsymbol{n}\mathcal{)}\simeq$ $\boldsymbol{n+1}$
\end{example}

\begin{example}
$\mathcal{I}\mathfrak{(}n\mathcal{)}\simeq$ $\mathcal{B}_{n}$, since every
subset of $n$ is an ideal.
\end{example}

If in a poset, $\mathcal{P}=(V;\leq),$ every pair of elements, $\left\{
x,y\right\}  $ has a least upper bound (greatest lower bound), it is denoted
$x\vee y$ ($x\wedge y$) and called the \textit{join }(\textit{meet}) of $x$
and $y.$ $\mathcal{L}=\left(  V;\vee,\wedge\right)  $ is then a
\textit{lattice}. Note that if $\mathcal{L}$ is finite it must have a least
element, $\bot$, and a greatest element, $\top$. Because of the roles they
play in the algebra of lattices, $\bot$ is often denoted as $0$, and $\top$ as
$1$. We prefer $\bot,\top$ because $0,1$ are already overloaded.

An element, $x\neq\bot$, in a lattice $\mathcal{L}=\left(  V;\vee
,\wedge\right)  $, is called \textit{join-irreducible} if%
\[
\nexists y,z<x\text{ such that }y\vee z=x.
\]

That is, $x$ has exactly one immediate predecessor.

\begin{example}
In the chain, $\boldsymbol{n=\left\{  0<1<...<n-1\right\}  ,}$ every element,
except $0$, is join-irreducible.
\end{example}

\begin{example}
In $\boldsymbol{2}^{V}$, the join-irreducible elements are exactly the
generators (singleton sets, elements of rank 1).
\end{example}

For a lattice, $\mathcal{L}=\left(  V;\vee,\wedge\right)  \mathcal{\ }$define
$\mathcal{J}\mathfrak{(}\mathcal{L)}$ to be $\left\{  x\in V:x\text{ is
join-irreducible in }\mathcal{L}\right\}  $, partially ordered by its induced
order in $\mathcal{L}$.

A lattice, $\mathcal{L}=\left(  V;\vee,\wedge\right)  $ is called
\textit{distributive} if $\forall x,y,z\in V$ it satisfies the

\begin{description}
\item[Distributive Laws] $x\wedge\left(  y\vee z\right)  =\left(  x\wedge
y\right)  \vee\left(  x\wedge z\right)  ,$ $x\vee\left(  y\wedge z\right)
=\left(  x\vee y\right)  \wedge\left(  x\vee z\right)  .$
\end{description}

\begin{example}
The Boolean lattice, $\mathcal{B}_{n}\simeq\left(  2^{n};\cup,\cap\right)  $
is distributive, in fact for any finite $\mathcal{P}$, $\mathcal{I}\left(
\mathcal{P}\right)  $ is closed under $\cup$ \& $\cap$ and inherits the
distributive laws from $\mathcal{B}_{n}$, where $n=\left\vert \mathcal{P}%
\right\vert $.
\end{example}

The following theorem is Birkhoff's fundamental result characterizing finite
distributive lattices. In \cite{D-P} it is Theorem 5.12.

\begin{theorem}
Let $\mathcal{L}=\left(  V;\vee,\wedge\right)  $ be a (finite) distributive
lattice. Then the map $\eta:\mathcal{L}\rightarrow\mathcal{I}\left(
\mathcal{J}\left(  \mathcal{L}\right)  \right)  $ defined by
\[
\eta\left(  x\right)  =\left\{  y\in\mathcal{J}\left(  \mathcal{L}\right)
:y\leq x\right\}
\]
is an isomorphism of $\mathcal{L}$ onto $\mathcal{I}\left(  \mathcal{J}\left(
\mathcal{L}\right)  \right)  .$
\end{theorem}

For a finite distributive lattice, $\mathcal{L}$, we call $\mathcal{I}%
\mathfrak{(}\mathcal{P)}$ (with $\mathcal{P}=\mathcal{\mathcal{J}\left(
\mathcal{L}\right)  }$) its \textit{Birkhoff representation}. Note that
$\mathcal{I}\left(  \mathcal{J}\left(  \mathcal{L}\right)  \right)  $ is a
sublattice of $\boldsymbol{2}^{V}$. Furthermore,

\begin{theorem}
$\mathcal{P\subseteq Q\Leftrightarrow I}\mathfrak{(}\mathcal{Q)\subseteq
I}\mathfrak{(}\mathcal{P)}$.
\end{theorem}

This is a special case of Theorem 5.19 of \cite{D-P}.

\section{General Derivation of StOp-order}

\subsection{The range of $\varphi:\mathcal{I}\mathfrak{(}%
\mathcal{P)\rightarrow I}\mathfrak{(}\mathcal{P)}$}

Now we come to the main result of this paper. Having stated Birkhoff's theorem
(Theorem 1 above) we can reveal the insight that lead to our result: The range
of $\varphi$ is a subposet of $\mathcal{I}\mathfrak{(}\mathcal{P)} $.
$\mathcal{I}\mathfrak{(}\mathcal{P)}$ in turn is a sublattice of the Boolean
lattice, $\boldsymbol{2}^{V}$ (closed under $\cup$ and $\cap$ and therefore
distributive). If $Range(\varphi)$ is representable as $\mathcal{I}%
\mathfrak{(}\mathcal{Q)}$ for some extension $\mathcal{Q}$ of $\mathcal{P}$,
then it would be closed under $\cup$ and $\cap$ and again, distributive. But
if we can just show that $Range(\varphi)$ is closed under $\cup$ and $\cap$
then it will be a sublattice of $\mathcal{I}\mathfrak{(}\mathcal{P)}$, must be
distributive and, by Theorems 1\&2, isomorphic to $\mathcal{I}\mathfrak{(}%
\mathcal{Q)}$ for some extension $\mathcal{Q}$ of $\mathcal{P}$.

First a preliminary result. We have observed that repeated composition of a
StOp with itself will produce an idempotent StOp with a smaller range so we
need only consider idempotent StOps.

\begin{lemma}
For an idempotent StOp $\varphi$, $Range(\varphi)=\left\{  T\in\mathcal{I}%
\mathfrak{(}\mathcal{P)}:\varphi\left(  T\right)  =T\right\}  $, the set of
fixpoints of $\varphi$.

\begin{proof}
By definition $\left\{  T\in\mathcal{I}\mathfrak{(}\mathcal{P)}:\varphi\left(
T\right)  =T\right\}  $ is a subset of $Range(\varphi)$. Conversely, $T\in
Range(\varphi)\Leftrightarrow\exists S\in\mathcal{I}\mathfrak{(}\mathcal{P)}$
such that $\varphi\left(  S\right)  =T$. But then
\begin{align*}
\varphi\left(  T\right)   &  =\varphi\left(  \varphi\left(  S\right)  \right)
\\
&  =\varphi^{2}\left(  S\right) \\
&  =\varphi\left(  S\right)  \text{, since }\varphi\text{ is idempotent,}\\
&  =T.
\end{align*}

\end{proof}

\begin{theorem}
For every idempotent Steiner operation, $\varphi:\mathcal{I}\mathfrak{(}%
\mathcal{P)}\rightarrow\mathcal{I}\mathfrak{(}\mathcal{P)}$, there exists a
unique partial order, $\mathcal{Q}$ on $V$ with $\mathcal{P\subseteq Q}$, such
that $Range(\varphi)=\mathcal{I}\mathfrak{(}\mathcal{Q)}$.
\end{theorem}

\begin{proof}
As remarked above, we need only show that the range of $\varphi$ is closed
under $\cup$ $\&$ $\cap$:
\begin{align*}
S,T\text{ }  &  \in\text{ }Range(\varphi)\\
&  \Rightarrow S,T\in\mathcal{I}\mathfrak{(}\mathcal{P)}\\
&  \Rightarrow S\cup T\in\mathcal{I}\mathfrak{(}\mathcal{P)}\text{, since
}\mathcal{I}\mathfrak{(}\mathcal{P)}\text{ is closed under }\cup.
\end{align*}
Also,%
\begin{align*}
\left(  S,T\subseteq S\cup T\right)   &  \Rightarrow\varphi\left(  S\right)
,\varphi\left(  T\right)  \subseteq\varphi\left(  S\cup T\right)  ,\text{ by
Axiom 3,}\\
&  \Rightarrow S\cup T\subseteq\varphi\left(  S\cup T\right)  \text{, since
}\varphi\left(  S\right)  =S\text{, }\varphi\left(  T\right)  =T\text{ }\\
&  \text{\& }\left(  S,T\subseteq W\Rightarrow S\cup T\subseteq W\right)
\text{.}%
\end{align*}
But by Axiom 1, $\left(  \left\vert \varphi\left(  S\cup T\right)  \right\vert
=\left\vert S\cup T\right\vert \right)  $ so $\varphi\left(  S\cup T\right)
=S\cup T$. Therefore $S\cup T\in Range(\varphi)$.

By duality, $S\cap T\in Range(\varphi)$.
\end{proof}
\end{lemma}

\section{\bigskip The Range of Ranges}

We wish to investigate the structure of all possible ranges of (finite)
Steiner operations. To this end we return to the strategy that led us to StOps
in the first place: We study morphisms for finite distributive lattices and
the resulting category. Morphisms for lattices are easy to define, they are
maps, $\varphi:\mathcal{L}\rightarrow\mathcal{M}$ that preserve the lattice
operations: $\forall x,y\in\mathcal{L}$,
\begin{align*}
\varphi\left(  x\wedge y\right)   &  =\varphi\left(  x\right)  \wedge
\varphi\left(  y\right)  ,\\
\varphi\left(  x\vee y\right)   &  =\varphi\left(  x\right)  \vee
\varphi\left(  y\right)  .
\end{align*}
If $\mathcal{L}$ is distributive and $\varphi$ is epi (onto), $\mathcal{M}$
must also be distributive. $\mathcal{L}$ is \textit{complete} means that
$\forall S\subseteq\mathcal{L},$ $%
{\displaystyle\bigvee}
S=x_{1}\vee x_{2}\vee...$ (for all $x_{i\in S}$) is defined and $%
{\displaystyle\bigwedge}
S$ is also defined (Definition 2.4(ii) of \cite{D-P}). Corollary 2.25 of
\cite{D-P} states that every finite lattice is complete. If $\mathcal{L}$ is a
complete lattice $%
{\displaystyle\bigvee}
\mathcal{L=}\top$ and $%
{\displaystyle\bigwedge}
\mathcal{L=\bot}$ so it has a top and bottom. We require that our lattice
morphisms preserve all existing meets and joins. And they must also preserve
top and bottom,
\begin{align*}
\varphi\left(  \top\right)   &  =\top,\\
\varphi\left(  \bot\right)   &  =\bot.
\end{align*}
In \cite{D-P} these are denoted \{0,1\}-homomorphisms, so we call them
\{$\bot$,$\top$\}-morphisms.

\subsection{The Category of Finite Distributive Lattices}

\begin{definition}
The category of distributive lattices with \{$\bot$,$\top$\}-morphisms will be
denoted DL.
\end{definition}

Note that $\boldsymbol{1}$ is the unique lattice$\boldsymbol{\ }$with $\perp$
$=\top$.

\begin{definition}
POSET is the category (see \cite{Mac}) whose objects are posets,
$\mathcal{P}=\left(  V,\leq\right)  $, and whose morphisms are monotone
(order-preserving) maps $\varphi:\mathcal{P\rightarrow Q}$. \textit{I.e.}
$\varphi$ is a function from $V_{\mathcal{P}}$ to $V_{\mathcal{Q}}$ such that
$x\leq_{\mathcal{P}}y$ $\Rightarrow\varphi\left(  x\right)  \leq_{\mathcal{Q}%
}\varphi\left(  y\right)  $.
\end{definition}

By the Connecting Lemma of \cite{D-P}, DL is (isomorphic to) a subcategory of
POSET but it is not full: Every lattice-morphism is order-preserving, but not
every order-preserving function between lattices preserves the lattice
operations (See Sections 2.16 to 2.19 of \cite{D-P}).

The restriction of POSET to finite posets is denoted POSET$_{F}$. Then an
extension of Birkhoff's theorem (Theorem 5.19 in \cite{D-P}) states that the
category DL$_{F}$ is isomorphic to POSET*$_{F}$, the dual of POSET$_{F}$.
Davey \& Priestley \cite{D-P} point out that the Birkhoff representation acts
a lot like the logarithm function of arithmetic, replacing large, apparently
complex structures in DL$_{F}$ by smaller ones in POSET*$_{F}$ but maintaining
their essential relationship. DL is a subcategory of POSET. The two categories
appear to differ considerably. However, Birkhoff's Theorem tells us that, when
restricted to finite posets and finite distributive lattices, they are
anti-isomorphic: $\mathcal{I}\mathfrak{:}$ POSET$_{F}$ $\rightarrow$ DL$_{F}$
and $\mathcal{J}:$ DL$_{F}\rightarrow$ POSET$_{F}$ are contravariant functors,
in fact $\mathcal{I}\left(  \mathcal{P}\right)  \simeq Hom_{POSET_{F}%
}(\mathcal{P}$, $\boldsymbol{2})$ by $I=\varphi^{-1}\left(  0\right)  $, is a
representable functor, and $\mathcal{J}$ is forgetful (of the join-reducible
elements and the lattice operations). $\mathcal{J}$ is essentially the inverse
of $\mathcal{I}$, \textit{i.e. }$\mathcal{J\circ I}=I_{POSET_{F}}$, the
identity functor on POSET$_{F}$ and $\mathcal{I\circ J}$ is naturally
isomorphic to $I_{\text{DL}_{F}}$, the identity functor on the category of
finite distributive lattices (See \cite{Gra}, Section II.1.3 for details).

\subsection{The Structure of DL$_{F}$}

We wish to determine universal constructions (limits and colimits) on the
category of distributive lattices, particularly the finite ones. As observed
by Davey \& Priestley \cite{D-P}, we need only study the category of posets.
In looking for a finite limit or colimit in POSET, a standard strategy is to
restrict the limit (or colimit) diagram to SET, in which all finite limits
exist, and try to show that the resulting limit in SET is actually a limit
(colimit) in POSET. This works for the basic limits, initial object, product
and equalizer. It also works for terminal object and coproduct, but not
coequalizer where a bit of tweeking is required (See \cite{Gou}, \cite{Awo}
(p. 126), and the Introduction of \cite{B-B-P} ). Anyway, POSET has all finite
limits and colimits.

\begin{example}
Let $\mathcal{P=}$ $\boldsymbol{1}$, $\mathcal{Q=}$ $\boldsymbol{3}$ and
define $\varphi_{1},\varphi_{2}:\mathcal{P\rightarrow Q}$ by $\varphi
_{1}\left(  0\right)  =0$ and $\varphi_{2}\left(  0\right)  =2$. Then in POSET
their coequalizer $\mathcal{C}$ has only one equivalence class ($\left\{
0,1,2\right\}  $) whereas in SET it has two ($\left\{  0,2\right\}  $ \&
$\left\{  1\right\}  $).
\end{example}

However, in many cases, particularly those that arise in applications,
coequalizers in POSET can be constructed as though they were in SET. Poset
morphisms that preserve the covering relation ($x\lessdot y\Rightarrow
\varphi\left(  x\right)  \lessdot\varphi\left(  y\right)  $) are equivalent to
digraph morphisms. DIGRAPH, the category of directed graphs, is a functor
category, DIGRAPH $\simeq$ FUNCT(D, SET) where D is the diagram category in
Figure 1 (appended). As a functor category, DIGRAPH inherits all limits from
its range, SET (Theorem 1, p. 115 of \cite{Mac}). If the limit is acyclic,
then it is the Hasse diagram for a partial order, the limit in POSET.

\subsubsection{Limits and Colimits in DL$_{F}$}

As in POSET, we can try to construct a given limit (or colimit) directly from
the corresponding limit (colimit) in SET. If that does not work we can
construct it as the image under $\mathcal{I}:$ POSET $_{F}\rightarrow$
DL$_{F}$ of the colimit (limit) of the image under the forgetful functor
$\mathcal{J}:$ DL$_{F}$ $\rightarrow$ POSET$_{F}$.

\begin{description}
\item[Initial Object] $\boldsymbol{2=}\left\{  0<1\right\}  $. Given any
finite lattice $\mathcal{L\in}$ DL, there is a unique lattice-morphism
$\varphi:\boldsymbol{2}\rightarrow\mathcal{L}$ defined by $\varphi\left\{
0\right\}  =$ $\perp$ \& $\varphi\left\{  1\right\}  =\top$. Note that
$\boldsymbol{2}$ is not the initial object of POSET ($\boldsymbol{0}$ is) but
$\boldsymbol{1}$ is the terminal object of POSET and $\boldsymbol{2=}%
\mathcal{I}\left(  \boldsymbol{1}\right)  .$

\item[Products] Given $\mathcal{L}$, $\mathcal{M\in}$ DL, $\mathcal{\ }$%
\[
\mathcal{L}\times\mathcal{M}=\left(  L_{\mathcal{L}}\times L_{\mathcal{M}%
};\wedge_{\mathcal{L}\times\mathcal{M}},\vee_{\mathcal{L}\times\mathcal{M}%
};\bot_{\mathcal{L}\times\mathcal{M}},\top_{\mathcal{L}\times\mathcal{M}%
}\right)  ,
\]
with lattice operations defined coordinatewise. \textit{I.e}.
\begin{align*}
\left(  x_{1},y_{1}\right)  \wedge_{\mathcal{L}\times\mathcal{M}}\left(
x_{2},y_{2}\right)   &  =\left(  x_{1}\wedge_{\mathcal{L}}y_{1},x_{2}%
\wedge_{\mathcal{M}}y_{2}\right)  ,\\
\left(  x_{1},y_{1}\right)  \vee_{\mathcal{L}\times\mathcal{M}}\left(
x_{2},y_{2}\right)   &  =\left(  x_{1}\vee_{\mathcal{L}}y_{1},x_{2}%
\vee_{\mathcal{M}}y_{2}\right)  ,\\
\bot_{\mathcal{L}\times\mathcal{M}}  &  =\left(  \bot_{\mathcal{L}}%
,\bot_{\mathcal{M}}\right)  \text{ and }\\
\top_{\mathcal{L}\times\mathcal{M}}  &  =\left(  \top_{\mathcal{L}}%
,\top_{\mathcal{M}}\right)  .
\end{align*}
The unique projection maps, $\pi_{1}:\mathcal{L\times M\rightarrow L}$ and
$\pi_{2}:\mathcal{L\times M\rightarrow M}$ are then DL-morphisms.

\item[Equalizers] Given parallel DL-morphisms $\varphi_{1}%
:\mathcal{L\rightarrow M}$ and $\varphi_{2}:\mathcal{L\rightarrow M}$,
$\mathcal{E}$ $=\left\{  x\in\mathcal{P}:\varphi_{1}\left(  x\right)
=\varphi_{2}\left(  x\right)  \right\}  $ with the embedding map,
$\varphi:\mathcal{E\rightarrow L}$ is their equalizer. It is trivially a DL-morphism.
\end{description}

Once again, by Mac Lane, DL$_{F}$ has all finite limits. Note that products
and equalizers are inherited from SET, but the initial object is not.

\begin{description}
\item[Terminal Object] The singleton lattice, $\boldsymbol{1}=\left\{
0\right\}  $. Given any $\mathcal{L\in}$ DL, there is a unique DL-morphism
$\boldsymbol{1}:\mathcal{L\rightarrow}\boldsymbol{1}$ defined by
$\boldsymbol{1}\left(  x\right)  =0$. Note that $\boldsymbol{1}$ is the
terminal object in POSET, but also $\boldsymbol{1=}\mathcal{I}\left(
\boldsymbol{0}\right)  $ and $\boldsymbol{0}$ is the initial object in POSET.

\item[Coproducts] Given $\mathcal{L}$, $\mathcal{M}\in$ DL, their coproduct,
in SET (and POSET) is the disjoint union of $\mathcal{L}$, $\mathcal{M}$. But
the disjoint union is not closed under meets and joins and there is no obvious
way to tweek it. However, if $\mathcal{L}$, $\mathcal{M}$ are finite and
distributive, Birkhoff's theorem gives $\mathcal{I}\left(  \mathcal{J}\left(
\mathcal{L}\right)  \times\mathcal{J}\left(  \mathcal{M}\right)  \right)  $ as
their coproduct in DL$_{F}$ . We denote it by $\mathcal{L}+\mathcal{M}$.

\item[Coequalizers] Given parallel DL-morphisms $\varphi_{1},\varphi
_{2}:\mathcal{L\rightarrow M}$, we have the same difficulty extending the
coequalizer from SET that we had in POSET. However, we have the equalizer,
\[
\mathcal{E}=\left\{  x\in\mathcal{J}\left(  \mathcal{M}\right)  :\mathcal{J}%
\left(  \varphi_{1}\right)  \left(  x\right)  =\mathcal{J}\left(  \varphi
_{2}\right)  \left(  x\right)  \right\}
\]
in POSET and Birkhoff's theorem guarantees that $\mathcal{I}\left(
\mathcal{E}\right)  $ will be their coequalizer in DL$_{F}$.
\end{description}

So again, by Mac Lane, DL$_{F}$ has all finite colimits, but this time only
the terminal object is inherited from SET.

\section{Natural Distributive Lattices}

We have shown that for a Steiner operation, $\varphi:\mathcal{I}%
\mathfrak{(}\mathcal{P)\rightarrow I}\mathfrak{(}\mathcal{P)}$, $Range(\varphi
)$ is a distributive lattice, $\mathcal{I}\mathfrak{(}\mathcal{Q)}$ for some
poset, $\mathcal{Q}$, that extends $\mathcal{P}$. We have argued that one
should only consider consistent StOps since they can be systematically
combined to get an even better StOp. Thus, given an isoperimetric problem on a
set, $V$, the possible ranges for StOps under consideration will be
sublattices of $2^{V}$ whose representing posets are suborders of a fixed
total order, $\tau$, of $V$. For purposes of study, we may take $V=n$ and the
total order to be the natural one, $\boldsymbol{n=}\left\{
0<1<...<n-1\right\}  $. This means that all the sublattices of $2^{n}$ will
contain $\mathcal{I}\mathfrak{(}\boldsymbol{n}\mathcal{)}=\boldsymbol{n+1}$.
Let us call these the \textit{natural distributive lattices (of order }$n$).
Ordered by $\subseteq$ they form a poset, $\mathcal{NDL}\left(  n\right)  $.
It is not hard to see that $\mathcal{NDL}\left(  n\right)  $ is a lattice with
$\mathcal{L\wedge M=L\cap M} $ and $\mathcal{L\vee M=}\overline{\mathcal{L\cup
M}}$, the closure of $\mathcal{L\cup M}$ under unions and intersections. Its
minimum element $\bot=$ $\boldsymbol{n+1}$ and its maximum element $\top=$
$2^{n}$. What else can we say about the structure of $\mathcal{NDL}\left(
n\right)  $?:

\begin{description}
\item[Question 1] Does it satisfy the Jordan-Dedekind chain condition?

\item[Question 2] If so, is it distributive or modular?
\end{description}

These look like challenging questions. Fortunately, following the recipe of
Davey \& Priestley in Section 4 of \cite{D-P}, we have an easy way to answer
them: A \textit{natural partial order (of order }$n$) is any suborder of
$\boldsymbol{n}$. Let $\mathcal{NPO}\left(  n\right)  $ be the set of all
natural partial orders (suborders of $n$) ordered by $\subseteq$. Birkhoff's
Theorem tells us that $\mathcal{NDL}\left(  n\right)  $ is isomorphic to
$\mathcal{NPO}^{\ast}\left(  n\right)  $, the dual of $\mathcal{NPO}\left(
n\right)  $. $\mathcal{NPO}\left(  n\right)  $ has already been investigated
by S. P. Avann \cite{Ava} and R. A. Dean \& G. Keller \cite{D-K}. We can
translate their findings into theorems about $\mathcal{NDL}\left(  n\right)  $:

\begin{description}
\item[Answer 1] \bigskip$\mathcal{NPO}\left(  n\right)  $ satisfies the
Jordan-Dedekind chain condition. The rank of $\mathcal{P}$ is
$r_{\mathcal{NPO}\left(  n\right)  }\left(  \mathcal{P}\right)  =\left\vert
<_{\mathcal{P}}\right\vert $, so $0\leq r_{\mathcal{NPO}\left(  n\right)
}\left(  \mathcal{P}\right)  $ $\mathcal{\leq}$ $\binom{n}{2}$ \cite{Ava}.
Therefore $\mathcal{NDL}\left(  n\right)  $ satisfies the Jordan-Dedekind
chain condition and its rank function is $r_{\mathcal{NDL}\left(  n\right)
}\left(  \mathcal{I}\left(  \mathcal{P}\right)  \right)  =r_{\mathcal{NPO}%
\left(  n\right)  }^{\ast}\left(  \mathcal{P}\right)  =\binom{n}%
{2}-r_{\mathcal{NPO}\left(  n\right)  }\left(  \mathcal{P}\right)  $.

\item[Answer 2] $\mathcal{NPO}\left(  n\right)  $ is not distributive or even
modular. However, it is lower semimodular, \textit{i.e. }$\forall\mathcal{P}
$,$\mathcal{Q\in NPO}\left(  n\right)  ,$
\[
r_{\mathcal{NPO}\left(  n\right)  }\left(  \mathcal{P}\right)
+r_{\mathcal{NPO}\left(  n\right)  }(\mathcal{Q})\leq r_{\mathcal{NPO}\left(
n\right)  }(\mathcal{P\wedge Q)}+r_{\mathcal{NPO}\left(  n\right)
}(\mathcal{P\vee Q})\text{ \cite{Ava}.}%
\]
Therefore $\mathcal{NDL}\left(  n\right)  $ is upper semimodular,
\textit{i.e.} the inequality above is reversed for $r_{\mathcal{NDL}\left(
n\right)  }$.
\end{description}

There are many more fascinating facts about natural partial orders (and thus
natural distributive lattices) in \cite{Ava} and \cite{D-K}. Also, the
sequence $\left\vert \mathcal{NPO}\left(  n\right)  \right\vert =\left\vert
\mathcal{NDL}\left(  n\right)  \right\vert $ is A006455 in The On-Line
Encyclopedia of Integer Sequences (OEIS \cite{OEIS} ). The second column of
the following table contains all known values:%

\[%
\begin{tabular}
[c]{c|c|c|r}%
$n$ & $\left\vert \mathcal{NPO}\left(  n\right)  \right\vert $ & $BPS\left(
n\right)  $ & $BPS\left(  n\right)  /\left\vert \mathcal{NPO}\left(  n\right)
\right\vert $\\\hline
\multicolumn{1}{r|}{$0$} & \multicolumn{1}{|r|}{$1$} &
\multicolumn{1}{|r|}{$0.0$} & $\allowbreak\allowbreak0.0$\\
\multicolumn{1}{r|}{$1$} & \multicolumn{1}{|r|}{$1$} &
\multicolumn{1}{|r|}{$15.\,\allowbreak179$} & $15.\,\allowbreak179$\\
\multicolumn{1}{r|}{$2$} & \multicolumn{1}{|r|}{$2$} &
\multicolumn{1}{|r|}{$51.\,\allowbreak055$} & $25.\,\allowbreak528$\\
\multicolumn{1}{r|}{$3$} & \multicolumn{1}{|r|}{$7$} &
\multicolumn{1}{|r|}{$\allowbreak182.\,\allowbreak14$} & $26.\,\allowbreak
02$\\
\multicolumn{1}{r|}{$4$} & \multicolumn{1}{|r|}{$40$} &
\multicolumn{1}{|r|}{$\allowbreak816.\,\allowbreak87$} & $\allowbreak
\allowbreak20.\,\allowbreak422$\\
\multicolumn{1}{r|}{$5$} & \multicolumn{1}{|r|}{$357$} &
\multicolumn{1}{|r|}{$4857.\,\allowbreak1$} & $13.\,\allowbreak605$\\
\multicolumn{1}{r|}{$6$} & \multicolumn{1}{|r|}{$4824$} &
\multicolumn{1}{|r|}{$39210.$} & $8.\,\allowbreak128\,1$\\
\multicolumn{1}{r|}{$7$} & \multicolumn{1}{|r|}{$96428$} &
\multicolumn{1}{|r|}{$4.\,\allowbreak352\,0\times10^{5}$} & $4.\,\allowbreak
513\,2$\\
\multicolumn{1}{r|}{$8$} & \multicolumn{1}{|r|}{$2800472$} &
\multicolumn{1}{|r|}{$6.\,\allowbreak691\,8\times10^{6}$} & $\allowbreak
\allowbreak2.\,\allowbreak389\,5$\\
\multicolumn{1}{r|}{$9$} & \multicolumn{1}{|r|}{$116473461$} &
\multicolumn{1}{|r|}{$1.\,\allowbreak432\,4\times10^{8}$} & $1.\,\allowbreak
229\,8$\\
\multicolumn{1}{r|}{$10$} & \multicolumn{1}{|r|}{$6855780268$} &
\multicolumn{1}{|r|}{$4.\,\allowbreak282\,8\times10^{9}$} & $0.624\,70$\\
\multicolumn{1}{r|}{$11$} & \multicolumn{1}{|r|}{$565505147444$} &
\multicolumn{1}{|r|}{$1.\,\allowbreak792\,8\times10^{11}$} & $0.317\,03$\\
\multicolumn{1}{r|}{$12$} & \multicolumn{1}{|r|}{$64824245807684$} &
\multicolumn{1}{|r|}{$\allowbreak1.\,\allowbreak052\,5\times10^{13}$} &
$\allowbreak0.162\,36$%
\end{tabular}
\]
$\ $.

Brightwell, Pr\"{o}mel \& Steger \cite{B-P-S} give a beautifully simple
formula,
\begin{align*}
BPS\left(  n\right)   &  =C_{n}n2^{\frac{n^{2}}{4}}\text{ with}\\
C_{n}  &  =\left\{
\begin{tabular}
[c]{rr}%
$12.7636300...$ & if $n$ is even\\
& \\
$12.7635965...$ & if $n$ is odd
\end{tabular}
\right.  \text{,}%
\end{align*}
and show that $BPS\left(  n\right)  /\left\vert \mathcal{NPO}\left(  n\right)
\right\vert \rightarrow1$ as $n\rightarrow\infty$. Since $C_{n}$ is the same
in the first 5 decimal places whether $n$ is even or odd, the difference does
not effect the values given in the table above. It is strange then, to see how
poor the approximation is for the known values (ratios are given in the fourth
column). This is not unprecedented however. It takes awhile for some of these
asymptotic sequences to settle down.

\subsection{What Finite Orders are StOp-Orders?}

We have shown that for any (finite) Steiner operation $\varphi$,
$Range(\varphi)$ is closed under $\cup$ $\&$ $\cap$ and is therefore a
distributive lattice. What other structure might $Range(\varphi)$ have? All
the StOp-orders in \cite{Har04} satisfy the Jordan-Dedekind chain condition.
Could that be a theorem? The following results answers this question in the negative.

\subsubsection{MWI Problems}

An interesting class of combinatorial isoperimetric problems is those for
which the boundary functional, $\omega:\mathcal{I}\left(  \mathcal{P}\right)
\rightarrow\mathbb{R}$, is additive:
\[
\omega\left(  S\right)  =%
{\displaystyle\sum\limits_{v\in S}}
\omega\left(  v\right)  \text{. }%
\]
Such an additive function is called a \textit{weight }and finding
\[
\min_{\substack{I\in\mathcal{I}\left(  \mathcal{P}\right)  \\\left\vert
I\right\vert =k}}\omega\left(  I\right)
\]
is called the \textit{minimum weight ideal }$\left(  \text{MWI}\right)  $
problem. Note that $\omega\left(  v\right)  $ may be negative as well as
positive, so
\[
\min_{\substack{I\in\mathcal{I}\left(  \mathcal{P}\right)  \\\left\vert
I\right\vert =k}}\omega\left(  I\right)  =-\max_{_{\substack{I\in
\mathcal{I}\left(  \mathcal{P}\right)  \\\left\vert I\right\vert =k}}}\left(
-\omega\left(  I\right)  \right)
\]
and minimizing or maximimizing are equivalent problems. If $\omega\left(
v\right)  <0$ for some $v\in V$ and $\min_{v\in V}\omega\left(  v\right)  =-C$
then $\omega^{+}\left(  v\right)  =\omega\left(  v\right)  +C\geq0$, and%
\[
\min_{\substack{I\in\mathcal{I}\left(  \mathcal{P}\right)  \\\left\vert
I\right\vert =k}}\omega^{+}\left(  I\right)  =\min_{\substack{I\in
\mathcal{I}\left(  \mathcal{P}\right)  \\\left\vert I\right\vert =k}%
}\omega\left(  I\right)  +kC
\]
so restricting $\omega$ to be positive makes no essential difference. The MWI
problem is trivial if $\mathcal{P=}\Delta_{V}$, the discrete order on $V $
($\mathcal{I}\left(  \Delta\right)  =2^{V}$): If we number the elements of
$V$, $\tau:V\rightarrow\left\{  1,2,...,n\right\}  $, one-to-one and onto, in
increasing order of their weight, $\tau\left(  u\right)  <\tau\left(
v\right)  \Rightarrow\omega\left(  u\right)  \leq\omega\left(  v\right)  $,
then $S_{m}=\left\{  \tau^{-1}\left(  1\right)  ,\tau^{-1}\left(  2\right)
,...,\tau^{-1}\left(  m\right)  \right\}  $ will be a solution of the MWI
problem. However, the general problem is NP-complete and many challenging
edge-isoperimetric and vertex-isoperimetric problems reduce to MWI problems
(see Section 6.2 of \cite{Har04}).

\subsubsection{Weight-Reductions}

As a special kind of the combinatorial isoperimetric problem, the defining
properties of Steiner operations apply to MWI problems. Stabilization and
compression are Steiner operations for MWI, but there is another systematic
family of Steiner operations that does not appear to apply to the EIP or VIP:
Suppose that $\mathcal{Q}$ is an extension of $\mathcal{P}=(V,\leq)$ and that
the weight function, $\omega$, for $\mathcal{P}$ is increasing on
$\mathcal{Q}$ ($\left(  u\leq_{\mathcal{Q}}v\right)  \Rightarrow\left(
\omega\left(  u\right)  \leq_{\mathcal{Q}}\omega\left(  v\right)  \right)  $).
Let $\tau:\mathcal{Q\rightarrow}\left[  n\right]  $ be any one-to-one \& onto
total extension of $\mathcal{Q}$, and $a$ any member of $V$. Then define
$\varphi_{a,\tau}:\mathcal{I}\left(  \mathcal{P}\right)  \rightarrow
\mathcal{I}\left(  \mathcal{P}\right)  $ by%
\[
\varphi_{a,\tau}\left(  I\right)  =I-v_{\max}+v_{\min},
\]
where $v_{\max}=\tau^{-1}\left(  \max\left\{  \tau\left(  v\right)
:a\leq_{\mathcal{Q}}v\in I\right\}  \right)  $ and $v_{\min}=\tau^{-1}\left(
\min\left\{  \tau\left(  u\right)  :a>_{\mathcal{Q}}u\notin I\right\}
\right)  $ if both sets are nonempty (if not, $\varphi_{a,\tau}\left(
I\right)  =I$). $I-v_{\max}\in\mathcal{I}\left(  \mathcal{P}\right)  $ because
$v_{\max}$ is maximal wrt $\mathcal{Q}$ and therefore wrt $\mathcal{P}$. By
the dual arguement, $(I-v_{\max})+v_{\min}\in\mathcal{I}\left(  \mathcal{P}%
\right)  $.

\begin{theorem}
$\varphi_{a,\tau}$ is a StOp for the MWI problem on $\left(  \mathcal{P}%
;w\right)  .$

\begin{proof}
We apply the definition of a StOp in Section 1.2:

\begin{enumerate}
\item $\left\vert \varphi_{a,\tau}\left(  I\right)  \right\vert =\left\vert
I\right\vert -1+1=\left\vert I\right\vert .$

\item $\omega\left(  \varphi_{a,\tau}\left(  I\right)  \right)  =\omega\left(
I\right)  -\omega\left(  v_{\max}\right)  +w\left(  v_{\min}\right)
\leq\omega\left(  I\right)  $ since $v_{\min}<_{\mathcal{Q}}a\leq
_{\mathcal{Q}}v_{\max}\Rightarrow\omega\left(  v_{\min}\right)  \leq
\omega\left(  v_{\max}\right)  $.

\item If $I\subseteq J$, then $\max\left\{  \tau\left(  v\right)
:a\leq_{\mathcal{Q}}v\in I\right\}  \leq\max\left\{  \tau\left(  v\right)
:a\leq_{\mathcal{Q}}v\in J\right\}  .$ If $=$ holds then $v_{\max}\left(
I\right)  =v_{\max}\left(  J\right)  $ and the same element is removed. from
$I,J$. If $<$ holds then an element not in $I$ will be removed from $J$. Also
$\min\left\{  \tau\left(  u\right)  :a>_{\mathcal{Q}}v\notin I\right\}
\leq\min\left\{  \tau\left(  u\right)  :a>_{\mathcal{Q}}u\notin J\right\}  $.
If $=$ holds then $v_{\min}\left(  I\right)  =v_{\min}\left(  J\right)  $ and
the same element is added to $I,J$. If $<$ holds then $v_{\min}\left(
I\right)  \in J$ already.

\item $\tau\left(  \varphi_{a,\tau}\left(  I\right)  \right)  =\tau\left(
I\right)  -\tau\left(  v_{\max}\right)  +\tau\left(  v_{\min}\right)  \leq
\tau\left(  I\right)  $, by the definition of $v_{\max}$, $v_{\min}$ and $=$
holds iff $\varphi_{a,\tau}\left(  I\right)  =I$.
\end{enumerate}
\end{proof}
\end{theorem}

Since $\varphi_{a,\tau}$ reduces (or at least does not increase) the weight of
an ideal, we call it a "\textit{reduction". }

\begin{theorem}
Every finite poset, $\mathcal{Q}$, is a StOp-order.

\begin{proof}
For a fixed $\tau$ the $\varphi_{a,\tau}$'s are consistent so the
superreduction, $\varphi_{\infty,\tau}$, defined by their cyclic composition
will be an idempotent StOp. $\varphi_{\infty,\tau}$ will determine a
StOp-order by Theorem 3 and it is easily seen that the StOp-order is
$\mathcal{Q}$.
\end{proof}
\end{theorem}

\begin{example}
Many of the StOps that we called "ad hoc" in \cite{Har04} are actually
reductions. Their definition seems superficial but the circumstances under
which they arise are still mysterious and they were useful in administering
the "coup de grace" after stabilization and compression had done the heavy
lifting. Anyway, those applications and Theorem 5 show that reductions are not
ad hoc but members of a rich systematic family of StOps.
\end{example}

\section{Conclusions \& Comments}

\subsection{Towards a Theory of StOp-orders}

Theorem 5 shows that in general StOp-orders have no additional structure.
However, many StOp-orders that occur in applications do have additional
structure. Some are distributive lattices themselves. Can such structure be
used to simplify their calculation? The Matsumoto-Verma theory of Bruhat
orders (the stabilization-orders derived from Coxeter groups), based on the
fact that Coxeter groups are generated by a relatively small subset of its
reflections (a basis) and that Bruhat orders have the Jordan-Dedekind chain
condition, is a great help in calculating Bruhat orders. Is there an extension
of those results to compression-orders or other families of Stop-orders?

\bigskip

\subsection{Do Continuous Steiner Operations Induce StOp-Orders?}

In 1966 \cite{Har66} the author labeled a combinatorial optimization problem
as "isoperimetric" because of its similarity with the classical isoperimetric
problem in the plane. The hope was that the analogy would guide intuition and
that techniques for classical (continuous) isoperimetric problems could be
extended to their combinatorial analogs. That hope has been fulfilled with the
theory of Steiner operations \cite{Har04}, applications of spectral theory
\cite{Chu} and abstract harmonic analysis \cite{C-H-P-P}. With Theorems 3 \& 5
it may now be possible for the combinatorial theory of Steiner operations to
repay something of its debt to classical analysis! Is there an analog of the
Birkhoff representation for the range of Steiner operations on continuous
measure spaces? For Steiner symmetrization the answer is, "yes, but the order
is not very interesting": The supersymmetrization of any bounded measureable
set in $\mathbb{R}^{n}$ is a sphere, centered at the origin, of the same
volume. Ordered by $\subseteq$, these spheres form a chain, isomorphic to
$\mathbb{R}_{+}$, which constitutes the symmetrization-order.

We began the search for a nontrivial StOp-order with Antonio Ros's survey of
classical isoperimetric problems \cite{Ros}. His paper was based on lectures
given at the Clay Mathematical Institute in 2001. In Section 1.6 Ros writes,
"The explicit description of the solutions of the isoperimetric problem in
flat 3-tori ($C_{1}^{3}$, the 3-fold product of unit circles) is one of the
nicest open problems in classical geometry". This is intriguing because it is
the $L_{2}$ analog of an $L_{1}$ problem solved by Bollobas \& Leader in 1991
(\cite{B-L91} or see \cite{Har04} Section 10.1). The Bollobas-Leader problem
is the continuous limit of the EIP on $\mathbb{Z}_{n}^{d}$ as $n\rightarrow
\infty$. They solved it in all dimensions $d$, even though it does not have
nested solutions for $d>1$, with a discontinuous variant of compression that
makes clever use of the convexity of the local solutions in dimension $d-1$.
Can the same strategy work for Ros's problem?

For $C_{1}$ ($d=1)$ the problem is trivial: The solutions are intervals of
length $v$, which may be nested. For $d>1$ we apply compression wrt this
1-dimensional solution and need only look at ideals in the product order of
$\left[  0,1\right]  ^{d}$ (note that there are just $d$ ways to factor
$C_{1}^{d}$ as a product $C_{1}\times C_{1}^{d-1}$). In addition we may apply
stabilization, the Steiner operation based on the reflective symmetries
induced by interchanging coordinates. The definition of stabilization for
continuous isoperimetric problems is the same as for combinatorial ones (see
Section 3.2.4 of \cite{Har04}) and is closely related to Hsiang symmetrization
(\cite{Ros}, Section 1.3). Stabilization may be made consistent with
compression and the resulting StOp-order is factorable as $\mathcal{L}\left(
d\right)  \times Stab\left(  Q_{d}\right)  $, where $\mathcal{L}\left(
d\right)  $is the standard simplex, $\left\{  0\leq x_{1}\leq x_{2}\leq...\leq
x_{d}\leq1\right\}  $ ordered coordinatewise, and $Stab\left(  Q_{d}\right)  $
is the stabilization-order of the graph of the $d$-cube, $Q_{d}$ (see Chapters
3 \& 4 of \cite{Har04}). The space of solutions is further limited by
regularity, which partitions it into components corresponding to the ideals of
$Stab\left(  Q_{d}\right)  $. When $d=3,$ $Stab\left(  Q_{3}\right)  $ has 10
ideals (Fig. 4.2 of \cite{Har04}) but the empty \& whole are trivial. Also,
throwing out those that are dual-complements of smaller ones (and therefore
redundant) we have the 5 on Ritor\'{e}'s list of candidates (\cite{Ros},
Section 1.5) when volume $v\leq1/2$. The same holds in any dimension but of
course it gets more complex (for $d=4$ (Fig. 4.3 of \cite{Har04}) there are 14 candidates).

There are further interaction between the local (variational) and global
(Steiner operational) conditions for an optimal surface: The regularity of the
surface implies it has a normal (directed outward) at every point$.$ In order
for the region enclosed to be an ideal in the StOp-order, that normal vector
must lie in the positive orthant. Also, where the surface meets the boundaries
of the cube, the normal vector must be orthogonal to the normal of the
bounding face. Where the surface intersects a hyperplane of symmetry of the
$d$-cube, $\left[  0,1\right]  ^{d}$, the dihedral angle between the tangent
plane and plane of symmetry must be acute (nonnegative inner product between
their normal vectors). The strongest local condition is that the mean
curvature of the surface must be constant.

In combinatorial isoperimetric problems (such as the edge-isoperimetric
problem on the graph of the d-cube) the notion of "nested solutions" is
fundamental. If one assumes that the problem might have nested solutions, then
starting with the empty set and adding in elements one at a time so as to
minimize the marginal boundary, with relatively little effort one has
candidates for solution sets of each cardinality. If those candidates
withstand scrutiny, then one has a powerful tool (compression) for proving
them optimal. The problems we are looking at now, however, are interesting
just because they do not have nested solutions. But they may have them in a
weaker sense. Suppose we start off with the ideal of volume 0 corresponding to
one of the ideals of Stab(Q\_d), thinking of it as an empty balloon. Pumping
air into the balloon will increase the volume so as to minimize the marginal
increase in area and should, intuitively, give a nested family of locally
optimal ideals. Conjecture 1 (Ritor\'{e}) of \cite{Ros07} affirms this
intuition in 3 dimensions. All we need to prove Ritor\'{e}'s conjecture is an
efficient analytic representation of all the locally optimal surfaces, but
evidently they do not exist for the larger two (Lawson's and Schwarz's
surfaces). One might also hope to adapt Bollobas \& Leader's discontinuous
compression arguement to prove the $L_{2}$ analog of their $L_{1}$ theorem.
However, we have not been able to do that either, and Ros thinks there might
be a counterexample\ in higher dimensions.

One of the notable controversies in mathematical history was Weierstrass's
challenge to Jakob Steiner's claim to have given a rigorous proof (the first)
of the classical isoperimetric theorem (\textit{circa} 1836, see \cite{Ber},
Section V.11., p. 295). Steiner showed that symmetrization of a planar set
(wrt a given line through its centroid) has the same area as the set and that
the length of its boundary is less (strictly) unless the set was already
symmetric \ (wrt the given line). Since the only planar set symmetric wrt
every line through its centroid is a circle, QED (Steiner claimed).
Weierstrass pointed out that Steiner was implicitly assuming that the
isoperimetric problem \textit{has} a solution and he still needed to prove it
for logical completeness. A proof of existence was finally published by
Schwartz in 1884. According to Berger \cite{Ber}, Blaschke's proof of
existence, based on a compactness argument, validated Steiner's intuition. We
are hopeful it can also prove that consistent Steiner operations generate
"pushouts". The basic symmetrizations and their "pushouts" are idempotent. Our
demonstration that the range of an idempotent Steiner operation is closed
under $\cup$ $\&$ $\cap$ does not invoke finiteness, so the ranges of the
basic- and super-symmetrizations will be (continuous) distributive lattices.
Our hope is that a variant of the Birkhoff-Priestley representation theory for
distributive lattices \cite{D-P} will produce its StOp-order. It seems that
there will have to be limitations on the closure of those lattices though,
like the countable unions of measure theory.

For a theory of continuous Steiner operations, the role played by Coxeter
groups (see \cite{Har04}, Chapter 5) should be taken by Lie groups. If the
action of a Lie group on a manifold is

\begin{enumerate}
\item Generated by reflections (order-two actions whose fixed submanifold
divides the manifold into two components),

\item Such that the stabilization it defines does not increase boundary,
\end{enumerate}

then the only subsets that need be considered in solving the isoperimetric
problem would be ideals in the stabilization-order (assuming that the Birkhoff
\& Priestley representation theories can produce a theoretical foundation for
such things).

\subsection{In Retrospect}

It might seem that Birkhoff's theorem was created to prove Theorem 3. However,
Birkhoff's theorem preceded Theorem 3 by at least 60 years. Also the essential
idea behind the proof of Birkhoff's theorem, the encoding of partial order
relations into the algebra of lattices (Theorem 2.8 (The Connecting Lemma) \&
Theorem 2.10 of \cite{D-P}), goes back another 60 years to Dedekind. The proof
of Theorem 3 is so simple (given Birkhoff's theorem and its extensions) yet
gives no insight into the interaction between StOps and the elements of the
underlying set (the join-irreducibles of $\mathcal{I}\left(  \mathcal{P}%
\right)  $). This lack of conceptual transparency indicates opportunity for
futher study.

\end{document}